\documentclass[11pt]{article}
\usepackage{amsthm}
\usepackage{latexsym,amssymb,amsmath}
\usepackage{graphicx,float,color,fancybox,shapepar,setspace,hyperref}
\usepackage{subfigure}
\usepackage{pgf,tikz}
\usetikzlibrary{arrows}
\newtheorem{thm}{Theorem}[section]
\newtheorem{conj}[thm]{Conjecture}

\newtheorem{prop}[thm]{Proposition}%[section]
 \def\dfn#1{{\sl #1}}
 \def\qed{\hfill\square}

\def\qed{ \hfill $\blacksquare$}
\def\pf{\medskip\noindent {\emph{Proof}.}~~}
\def\less{\setminus}

\newcounter{counter}
\begin{document}

\title{A note on Gallai-Ramsey number of even wheels}

\author{ Zi-Xia Song$^{a,}$\thanks{Supported by the National   Science  Foundation under Grant No. DMS-1854903.  
     E-mail address:  Zixia.Song@ucf.edu.}, Bing Wei$^b$, Fangfang Zhang$^{c,a,}$\thanks{Corresponding author.  The work is done while the third author was studying at the University of Central Florida as a visiting student, supported by the Chinese Scholarship Council.     E-mail address:   fangfangzh@smail.nju.edu.cn}\,   and Qinghong Zhao$^b$   \\
$^a${\small Department  of Mathematics, University of Central Florida, Orlando, FL 32816, USA} \\
$^b${\small Department  of Mathematics,    University of Mississippi, Oxford,  MS 38677,  USA} \\
$^c${\small Department of Mathematics,   Nanjing University,    Nanjing 210093, P. R. China}\\
 }

\maketitle

\begin{abstract}
A Gallai coloring of a complete graph is an edge-coloring such that no triangle   has  all its edges colored differently.  A Gallai $k$-coloring is a Gallai coloring that uses $k$ colors.   Given a graph $H$ and an integer $k\geq 1$,   the  Gallai-Ramsey number $GR_k(H)$ of $H$ is the least positive  integer $N$ such that every Gallai $k$-coloring of the complete graph $K_N$ contains a monochromatic copy of $H$.   Let $W_n $ denote a wheel  on  $n+1$ vertices. In this note, we study Gallai-Ramsey number of $W_{2n}$ and  completely determine  the exact value of   $GR_k(W_4)$ for all $k\ge2$. \\
 

\end{abstract}

\section{Introduction} 

All graphs in this paper   are finite, simple and undirected.    Given a graph $G$ and a set $S\subseteq V(G)$,  we use   $|G|$    to denote  the  number
of vertices    of $G$, and  $G[S]$ to denote the  subgraph of $G$ obtained from $G$ by deleting all vertices in $V(G)\less S$.    For  two disjoint sets $A, B\subseteq V(G)$,  $A$ is \dfn{complete} to $B$ in $G$  if every vertex in $A$ is adjacent to all vertices in  $B$.    We use $K_n, C_n, P_n$   to denote the
complete graph, cycle, and path on $n$ vertices, respectively; and   $W_n $ to denote a  wheel  on  $n+1$ vertices. 
For any positive integer $k$, we write  $[k]$ for the set $\{1, \ldots, k\}$. We use the convention   ``$S:=$'' to mean that $S$ is defined to be the right-hand side of the relation.\medskip

 Given an integer $k \ge 1$ and graphs $H_1,  \ldots, H_k$, the classical Ramsey number $R(H_1,   \ldots, H_k)$   is  the least    integer $N$ such that every $k$-coloring of  the edges of  $K_N$  contains  a monochromatic copy of  $H_i$ in color $i$ for some $i \in [k]$. When $H = H_1 = \dots = H_k$, we simply write    $R_k(H)$.   Ramsey numbers are notoriously difficult to compute in general.  Very little is known for Ramsey numbers of wheels and  the exact value of  $R_2(W_n)$ is only known for    $n\in\{3,4,5\}$ (see  \cite{W5,W3, R2}).    In this paper, we  study Ramsey number  of wheels in Gallai colorings, where a \dfn{Gallai coloring} is a coloring of the edges of a complete graph without rainbow triangles (that is, a triangle with all its edges colored differently). Gallai colorings naturally arise in several areas including: information theory~\cite{KG}; the study of partially ordered sets, as in Gallai's original paper~\cite{Gallai} (his result   was restated in \cite{Gy} in the terminology of graphs); and the study of perfect graphs~\cite{CEL}.  
  More information on this topic  can be found in~\cite{FGP, FMO}.  \medskip
 
A \dfn{Gallai $k$-coloring} is a Gallai coloring that uses at most $k$ colors.  We use $(G, \tau)$ to denote a Gallai $k$-colored complete graph if $G$ is a complete graph and   $\tau: E(G)\rightarrow [k]$ is a Gallai $k$-coloring.   Given an integer $k \ge 1$ and graphs $H_1,  \ldots, H_k$, the   \dfn{Gallai-Ramsey number} $GR(H_1,  \ldots, H_k)$ is the least integer $N$ such that  every   $(K_N, \tau)$   contains a monochromatic copy of $H_i$ in color $i$ for some $i \in [k]$. When $H = H_1 = \dots = H_k$, we simply write $GR_k(H)$.    Clearly, $GR_k(H) \leq R_k(H)$ for all $k\ge1$ and $GR(H_1, H_2) = R(H_1, H_2)$.    In 2010, 
Gy\'{a}rf\'{a}s,   S\'{a}rk\"{o}zy,  Seb\H{o} and   Selkow~\cite{exponential} proved   the general behavior of $GR_k(H)$.

\begin{thm} [\cite{exponential}]
Let $H$ be a fixed graph  with no isolated vertices 
 and let $k\ge1$ be an integer. Then
$GR_k(H) $ is exponential in $k$ if  $H$ is not bipartite,    linear in $k$ if $H$ is bipartite but  not a star, and constant (does not depend on $k$) when $H$ is a star.			
\end{thm}

It turns out that for some graphs $H$ (e.g., when $H=C_3$),  $GR_k( H)$ behaves nicely, while the order of magnitude  of $R_k(H)$ seems hopelessly difficult to determine.  It is worth noting that  finding exact values of $GR_k (H)$ is  far from trivial, even when $|H|$ is small.
 The following  structural result of Gallai~\cite{Gallai} is crucial in determining the exact value  of Gallai-Ramsey number of a graph $H$.

\begin{thm}[\cite{Gallai}]\label{Gallai}
	For any Gallai coloring $\tau$ of a complete graph $G$ with $|G| \ge 2$, $V(G)$ can be partitioned into nonempty sets  $V_1, V_2, \dots, V_p$ with $p>1$ so that    at most two colors are used on the edges in $E(G)\less (E(V_1)\cup \cdots\cup  E(V_p))$ and only one color is used on the edges between any fixed pair $(V_i, V_j)$ under $\tau$, where $E(V_i)$ denotes the set of edges in $G[V_i]$ for all $i\in [p]$. 
\end{thm}

The partition given in Theorem~\ref{Gallai} is  a \dfn{Gallai partition} of  the complete graph $G$ under $\tau$.  Given a Gallai partition $V_1,  \dots, V_p$ of the complete graph $G$ under  $\tau$, let $v_i\in V_i$ for all $i\in[p]$ and let $\mathcal{R}:=G[\{v_1,  \dots, v_p\}]$. Then $\mathcal{R}$ is  the \dfn{reduced graph} of $G$ corresponding to the given Gallai partition under  $\tau$. 
By Theorem~\ref{Gallai},  all edges in $\mathcal{R}$ are colored by at most two colors under  $\tau$.  One can see that any monochromatic copy of $H$ in $\mathcal{R}$  will result in a monochromatic copy of $H$ in $G$ under  $\tau$. It is not  surprising  that  Gallai-Ramsey number $GR_k( H)$ is closely  related to  the classical Ramsey number $R_2( H)$.  The following is a conjecture of   Fox,  Grinshpun and  Pach~\cite{FGP}   on $GR_k( H)$ when $H$ is a complete graph. 

\begin{conj}[\cite{FGP}]\label{Fox} For all  $k\ge1$ and $t\ge3$,
\[
GR_k( K_t) = \begin{cases}
			(R_2(K_t)-1)^{k/2} + 1 & \text{if } k \text{ is even} \\
			(t-1)  (R_2(K_t)-1)^{(k-1)/2} + 1 & \text{if } k \text{ is odd.}
			\end{cases}
\]
\end{conj}

The first case $t=3$ of Conjecture~\ref{Fox}  follows directly from a result of  Chung and Graham~\cite{chgr} from  1983, and a simpler proof  can be found in~\cite{exponential}. 
 Recently, the case $t=4$ of Conjecture~\ref{Fox} was proved in \cite{k4}. 
\medskip

   In this note, we study Gallai-Ramsey number  of even wheels.  Recall that  $(G, \tau)$   denotes a Gallai $k$-colored complete graph if $G$ is a complete graph and   $\tau: E(G)\rightarrow [k]$ is a Gallai $k$-coloring.  For any $A, B\subseteq V(G)$, if all the edges between $A$ and $B$ are colored the same color, say blue, under $\tau$,  we say that $A$ is \dfn{blue-complete} to $B$ in $(G, \tau)$. We simply say $a$ is  \dfn{blue-complete} to $B$ in $(G, \tau)$ when $A=\{a\}$. We first prove  the following   general lower bound for $GR_k(W_{2n})$ for all $k\ge2$ and $n\ge2$.

 \begin{prop} \label{p}
For   all  $n \ge 2$ and  $k\ge2$, 
$$GR_k(W_{2n}) \ge \begin{cases}
			(R_2(W_{2n})-1)\cdot 5^{(k-2)/2} + 1 & \text{if } k \text{ is even} \\
			2(R_2(W_{2n})-1) \cdot  5^{(k-3)/2} + 1 & \text{if } k \text{ is odd.}
			\end{cases}$$

\end{prop}

\pf   Let $n\ge2$ be given. For   all   $k\ge2$, let 
$$f(k):=   \begin{cases}
			(R_2(W_{2n})-1)\cdot 5^{(k-2)/2}   & \text{if } k \text{ is even} \\
			2(R_2(W_{2n})-1) \cdot  5^{(k-3)/2}   & \text{if } k \text{ is odd.}
			\end{cases}$$
For all $k\ge 2$, we next construct   
 $(G_k, \tau_k)$ recursively as follows,  where $G_k$  is a  complete graph on $f(k)$ vertices and  $\tau_k: E(G_k)\rightarrow [k]$ is   a Gallai $k$-coloring  such that $(G_k, \tau_k)$ has no monochromatic copy of $W_{2n}$.  \medskip
  
  For $k=2$, let $m:=f(2)=R_2(W_{2n})-1$ and $G_2:=K_{m}$.  By the definition of $R_2(W_{2n})$, there exists $\tau_2:E(G_2)\rightarrow \{1,2\}$ such that  $G_2$ has no monochromatic copy of $W_{2n}$ under $\tau_2$. Then $(G_2, \tau_2)$  is a  Gallai $2$-colored complete graph on $f(2)$ vertices  with      no monochromatic copy of $W_{2n}$, as desired. When  $k\ge3$ is even, let   $(G_{k-2}, \tau_{k-2})$ be the construction on $f(k-2)$ vertices  with colors in $[k-2]$ and let $(H, c):=(K_5, c)$  be such that   edges of $H$ are colored by colors $k$ and $k-1$ with  no monochromatic copy of $K_3$. Let $(G_k, \tau_k)$ be obtained from $(H, c)$ by replacing each vertex of $H$ with a copy of $(G_{k-2}, \tau_{k-2})$  such that for every $uv\in E(H)$, all the edges between the corresponding copies of $(G_{k-2}, \tau_{k-2})$ for $u$ and $v$ are colored by the color $c(uv)$.  Since $(G_{k-2}, \tau_{k-2})$ is a Gallai $(k-2)$-colored complete graph  with no monochromatic copy of $W_{2n}$, we see that $\tau_k$ is indeed a Gallai $k$-coloring and $(G_k, \tau_k)$ has no monochromatic copy of $W_{2n}$. When  $k\ge3$ is odd, let $(G_{k-1}, \tau_{k-1})$ be the construction on $f(k-1)$ vertices  with colors in $[k-1]$ and let $(G_k, \tau_k)$ be the join of two disjoint copies of $(G_{k-1}, \tau_{k-1})$, with all the new edges   colored by the color $k$.   In both cases, $(G_k, \tau_k)$ has no rainbow triangle and no monochromatic copy of $W_{2n}$.  Hence,  $GR_k(W_{2n}) \ge f(k)+1$, as desired. \qed\\

We will make use of the following  result of Hendry~\cite{R2} on the exact value of  $R_2(W_4)$.
\begin{thm}[\cite{R2}]\label{W4}
  $R_2(W_4)=15$.  
  \end{thm}\medskip
  
 The main result in this note is Theorem~\ref{thm} below, which establishes the exact value of $GR_k(W_4)$ for all $k\ge2$. \medskip
 
  \begin{thm}\label{thm} 
  
For all $k\ge 2$,
\[
GR_k(W_4)= \left\{
\begin{array}{cc}
14\cdot 5^{\frac{k-2}{2}}+1,& \text{ if }  k  \text{ is even}\\
28\cdot 5^{\frac{k-3}{2}}+1,& \text{ if }  k  \text{ is odd}.
\end{array}
\right.
\]
\end{thm}

 With the support of Theorem~\ref{thm}, it seems that the lower bound given in Proposition~\ref{p} is also the desired upper bound for $GR_k(W_{2n})$. We propose the following conjecture.
  
  \begin{conj}\label{CONJ} For all $k\ge 2$ and $n\ge 2$,
\[
GR_k(W_{2n})= \left\{
\begin{array}{cc}
(R_2(W_{2n})-1)\cdot 5^{\frac{k-2}{2}}+1,& \text{if}\ k\text{ is even}\\
2(R_2(W_{2n})-1)\cdot 5^{\frac{k-3}{2}}+1,&  \text{if}\ k\text{ is odd}.
\end{array}
\right.
\]
\end{conj}

 It is worth noting that  very recently,    Gallai-Ramsey number  of some other graphs with at most five vertices were studied by  Li and Wang \cite{five}. More recent work on   Gallai-Ramsey numbers of  cycles can be found in  \cite{C9C11,C13C15, DylanSong,  chen,  C10C12, C6C8, oddcycle}.

\section{Proof of Theorem \ref{thm}}

Let $f(1):=4$ and for each integer $s\ge2$,  let  \[
f(s):= \left\{
\begin{array}{cc}
14\cdot 5^{\frac{s-2}{2}},& \text{ if }  s  \text{ is even}\\
28\cdot 5^{\frac{s-3}{2}},& \text{ if } s  \text{ is odd}.
\end{array}
\right.
\]
 Then for all $s\ge3$, 
\[\ \ \ \ \ \ \ \ \ \ \ \ \
f(s)\ge  \left\{
\begin{array}{cc}
2f(s-1)\\
5f(s-2).
\end{array}
\right.
\tag{$\ast$}
\]

By Proposition \ref{p}, 
it suffices to show that for all  $k\ge 2$, 
 $GR_k(W_4)\le f(k)+1$.    
By Theorem~\ref{W4},  $ R_2(W_4)=15$ and so 
$GR_2(W_4)=R_2(W_4) =f(2)+1$.  We may assume that $k\ge 3$.   Let $G:=K_{f(k)+1}$ and let $\tau: E(G)\xrightarrow{} [k]$ be any Gallai $k$-coloring of $G$. Suppose that $(G, \tau)$ does not contain any monochromatic copy of $W_4$. Then $\tau$ is \dfn{bad}. Among all complete graphs on $f(k)+1$ vertices with a bad Gallai $k$-coloring, we choose $(G, \tau)$ with $k$ minimum. By the minimality of $k$, we have $GR_\ell(W_4)\le f(\ell)+1$ for all $  \ell\in\{2, \ldots,  k-1\}$. We next prove several claims.     \\

\noindent {\bf Claim\refstepcounter{counter}\label{H}  \arabic{counter}.}  
Let   $A\subseteq V(G)$ and let $\tau_\ell$ be a Gallai $\ell$-coloring of $G[A]$ with no monochromatic copy of $W_4$, where  $ 3\le \ell \le k$. Let $i$ be a color used by $\tau_\ell$.  If $(G[A], \tau_\ell)$ has no monochromatic copy of   $P_3$ in color $i$, then  $|A|\le f(\ell-1)$.

\pf Suppose    $|A|\ge f(\ell-1)+1$.  By the minimality of $k$, $GR_{\ell-1}(W_4)\le f(\ell-1)+1$. Then $|A|\ge GR_{\ell-1}(W_4)$.  We may assume that the color $i$ is blue. Since $\ell\ge3$, let red  be another color in $[\ell]\less \{i\}$, and let  $\tau^*$ be  obtained from $\tau_\ell$ by replacing colors red and blue   with a new color, say $\beta\notin [\ell]$.  Then $\tau^*$  is  a Gallai $(\ell-1)$-coloring of $G$ and so $(G[A], \tau^*)$ contains a monochromatic copy of $W:=W_4$ because $|A|\ge GR_{\ell-1}(W_4)$. Since   $(G[A], \tau_\ell)$ contains no monochromatic $W_4$, it follows that  $ W$ must be  a monochromatic  $ W_4 $ in color $\beta$. For the remainder of the  proof of this claim,  we use  $W_4=(u_1,  u_2, u_3, u_4; u_0)$ to denote a wheel  on  $5$ vertices, where  $W_4\less u_0$ is a cycle with vertices $u_1,  u_2, u_3, u_4$ in order, and $u_0$ is complete to $\{u_1,  u_2, u_3, u_4\}$.  Let $W= (u_1,  u_2, u_3, u_4; u_0)$.   Then $W$ has at most two independent  blue edges under  $\tau_\ell$, because    $(G[A], \tau_\ell)$ has no blue $P_3$.    Next, we show that  $W[\{u_1,  u_2, u_3, u_4\}]$ is  a red $C_4$ under  $\tau_\ell$. Without loss of generality, suppose that     $u_1u_2$ is blue. Since $(G[A], \tau_\ell)$ has no blue $P_3$, we see that $u_1u_4, u_1u_0, u_2u_3, u_2u_0$ must be red. We further observe that $u_1u_3$ and $u_2u_4$  must be red because $(G[A], \tau_\ell)$ has no rainbow triangle. It follows that  $\{u_1, u_2\}$ is red-complete to $\{u_3, u_4, u_0\}$ in $G$  under  $\tau_\ell$.  Note that  one of $u_0u_3$ and $u_0u_4$ is red since $(G[A], \tau_\ell)$ has no blue $P_3$.  We may further assume that $u_0u_3$ is red under  $\tau_\ell$.  But then we obtain  a red   $W_4=(u_1, u_4, u_2, u_3; u_0)$ under  $\tau_\ell$ when $u_0u_4$ is red,  and a red $W_4=(u_0, u_2, u_4, u_1; u_3)$ under  $\tau_\ell$ when $u_0u_4$ is blue, contrary to the fact that $(G,\tau_\ell)$ has no monochromatic copy of $W_4$.   This proves that $W[\{u_1,  u_2, u_3, u_4\}]$  is  a red $C_4$ under  $\tau_\ell$.  Since $(G[A], \tau_\ell)$ has no blue $P_3$, we see that  exactly   one of the edges  $u_0u_1, u_0u_2, u_0u_3, u_0u_4$, say $u_0u_1$ is blue under  $\tau_\ell$. But then we obtain  a red $W_4 =(u_1, u_2, u_0, u_4; u_3)$ under  $\tau_\ell$, contrary to the fact that $(G,\tau_\ell)$ has no monochromatic copy of $W_4$.   \qed\\

 \noindent {\bf Claim\refstepcounter{counter}\label{HH}  \arabic{counter}.}  
 Let   $A\subseteq V(G)$ and let  $i, j\in [k]$  be two distinct colors.  If $(G[A], \tau)$ has  no  monochromatic copy of  $P_3$ in color $i$  or    in color $j$,   then $|A|\le f(k-2)$.

\pf We may assume that  color $i$ is  red and color $j$ is blue. Then $(G[A], \tau)$ has  neither red nor blue    copy of  $P_3$.  Let  $\tau^*$ be  obtained from $\tau$ by replacing colors red  and blue with a new color, say $\beta\notin [k]$.   Then $\tau^*$  is   a Gallai $(k-1)$-coloring of $G$. We claim that $(G[A], \tau^*)$ has no monochromatic copy of $P_3$ in color $\beta$. Suppose not. Let $x_1, x_2, x_3\in A$ be such that $x_1x_2, x_2x_3$ are colored by color $\beta$ under $\tau^*$. Since $(G[A], \tau)$ has neither red nor blue $P_3$, we may assume that $x_1x_2$ is colored blue and $x_2x_3$ is colored red under  $\tau$. Then $x_1x_3$ is colored neither red nor blue under  $\tau$. But then $(G[\{x_1, x_2, x_3\}], \tau)$ is a rainbow triangle, a contradiction. Thus $(G[A], \tau^*)$ has no monochromatic copy of $P_3$ in color $\beta$, as claimed. Then $(G[A],\tau^*)$ has no monochromatic copy of  $W_4$.  By Claim~\ref{H} applied to $A$, $\tau^*$ and color $\beta$, we have $|A|\le f(k-2)$. \qed\\

  \noindent {\bf Claim\refstepcounter{counter}\label{Two}  \arabic{counter}.}  Let $A\subseteq V(G)$ and let $x, y\in V(G)\less A$ be distinct.   If all the edges between $\{x,y\}$  and  $A$ are colored the same color, say color $i\in[k]$, under  $\tau$,  then $(G[A], \tau)$ contains no monochromatic copy of  $P_3$ in color $i$. 

\pf We may assume that the color $i$ is blue. Suppose  $(G[A],\tau)$ contains a blue copy of $P_3$   with vertices, say  $v_1,v_2,v_3$ in order. Then  we obtain a blue $W_4=(v_1,x,v_3,y;  v_2)$, which is a  contradiction.\qed\\

Let $x_1,x_2,\ldots,x_m\in V(G)$ be a maximum sequence of vertices chosen as follows: for each 
$j\in [m]$, all edges between $x_j$ and $V(G)\less\{x_1,x_2,\ldots,x_j\}$ are colored by the same color under  $\tau$. Let $X:=\{x_1,x_2,\ldots, x_m\}$. Note that $X$ is possibly empty. For each $x_j\in X$, let $\tau(x_j)$ be the unique color on the edges between $x_j$ and $V(G)\less\{x_1,x_2,\ldots,x_j\}$.\\

 \noindent {\bf Claim\refstepcounter{counter}\label{X}  \arabic{counter}.}  
$\tau(x_i)\neq \tau(x_j)$ for all $i,j\in [m]$ with $i\neq j$.

\pf Suppose there exist $i,j\in [m]$ with $i\neq j$ such that $\tau(x_i)= \tau(x_j)$. We may assume that   $x_j$ is the first vertex in the sequence $x_1,\ldots,x_m$ such that $\tau(x_j)=\tau(x_i)$ for some $i\in [m]$ with $i<j$. We may further assume that the color $\tau(x_i)$ is blue. By the pigeonhole principle, $j\le k+1$. Let $A:=V(G)\less \{x_1,x_2,\ldots,x_j\}$. Then all the edges between $\{x_i, x_j\}$ and $ A$ are colored blue under $\tau$. By Claim~\ref{Two},   $(G[A],\tau)$ has no blue $P_3$.  By Claim \ref{H}, $|A|\le f(k-1)$. Then $|G|\le f(k-1)+k+1<f(k)+1$, which is impossible.  \qed\\

By Claim \ref{X}, $|X|\le k$. Note that  $G\less X$ has no monochromatic copy of $W_4$ under  $\tau$. Consider a  Gallai-partition of $G\less X$, as given in Theorem~\ref{Gallai},  with parts $V_1, V_2, \dots, V_{p}$, where $p\ge2$. We may assume that  $|V_1|\ge|V_2|\ge \ldots \ge |V_{p}|$.
Let $\mathcal{R}$ be  the reduced graph of $G\less X$ with vertices $a_1,\ldots,a_{p}$. 
By Theorem \ref{Gallai}, we may assume that every edge of $\mathcal{R}$ is colored red or blue. 
Note that any monochromatic copy of $W_4$ in $\mathcal{R}$ would yield a monochromatic copy of $W_4$ in $G\less X$. Thus $\mathcal{R}$ has neither red nor blue $W_4$. Since  $R_2(W_4)=15$, we see that $p\le 14$. Then  $|V_1|\ge 2$ because $|G\less X|\ge 26$  for all  $k\ge 3$. Let 
\[
\begin{split}
B &:= \{v\in \cup_{i=2}^p V_i \mid  v~\text{is blue-complete to}~V_1  \text{ under  }  \tau\} \text{ and} \\
R &:= \{v\in \cup_{i=2}^p V_i \mid  v~\text{is red-complete to}~V_1  \text{ under  }  \tau \}.
\end{split}
\]
 Then $|B|+|R|+|V_1|=|G|-|X|$.  We may assume that $|B|\le |R|$. Suppose $|R|=1$. Then $|B|\le 1$. We may assume that the only edge between $R$ and $B$ is red if   $|B|=1$. Then  all the edges between  the vertex in $R$ and $V(G)\less (X\cup R)$ are colored red under $\tau$, contrary to the maximality of $m$ when choosing   $x_1, \ldots, x_m$. Thus    $|R|\ge 2$.  By Claim~\ref{Two}, neither $(G[V_1], \tau)$ nor  $(G[R], \tau)$ contains a     red copy of $P_3$. Note that   $(G[V_1\cup R], \tau)$ has a red copy of $C_4$. Thus no vertex in $X$ is red-complete to $V(G)\less X$ under  $\tau$.  \\

 \noindent {\bf Claim\refstepcounter{counter}\label{V1}  \arabic{counter}.}  
$(G[B], \tau)$ contains a  red copy of $P_3$. Furthermore,   $|V_1\cup X|\le  f(k-2)$.

\pf Suppose $(G[B], \tau)$ has no  red $P_3$. Since   no vertex in $X$ is red-complete to $V(G)\less X$ under  $\tau$  and $(G,\tau)$ has no rainbow triangle, we see that $(G[X],\tau)$ has no red edge,  and no edge between $X$ and $V(G)\less X$ is colored red under $\tau$. Then    $(G[B\cup V_1\cup X], \tau)$ has no red $P_3$, because  $(G[V_1], \tau)$ has no red $P_3$.  Note that  $(G[R], \tau)$ has no red $P_3$. By Claim \ref{H}, $|R|\le f(k-1)$ and $|B\cup V_1\cup X|\le f(k-1)$. By ($\ast$), $2f(k-1)\le f(k)$ for all $k\ge3$.  But then $|G|=|R|+|B\cup V_1\cup X|\le 2f(k-1)\le f(k)$, which contradicts the fact that $|G|= f(k)+1$. Thus  $(G[B], \tau)$ contains a  red copy of $P_3$ and so $|B|\ge3$. If  some vertex  $x\in X$ is blue-complete to $V(G)\less X$ under  $\tau$, then $ (G,\tau)$ contains  a blue $W_4=(v_1, b_1, v_2, b_2; x)$, where $v_1, v_2\in V_1$, and $b_1, b_2\in B$. Thus  no vertex in $X$ is blue-complete to $V(G)\less X$ under  $\tau$. It follows that    $(G[X],\tau)$ has no blue edge,  and no edge between $X$ and $V(G)\less X$ is colored blue  under $\tau$.  By Claim~\ref{Two},   $(G[V_1], \tau)$  has no   blue $P_3$.  Thus $(G[V_1\cup X],\tau)$ contains neither red   nor blue $P_3$. By Claim \ref{HH}, $|V_1\cup X|\le  f(k-2)$, as desired. \qed\\

 \noindent {\bf Claim\refstepcounter{counter}\label{V3}  \arabic{counter}.}  
$|V_3|\le 1$.

\pf Suppose  $|V_1|\ge|V_2|\ge |V_3|\ge 2$. Recall that $a_1,a_2,a_3$ are the corresponding vertices of $V_1, V_2, V_3$ in the reduced graph $\mathcal{R}$. If $\mathcal{R}[\{a_1,a_2,a_3\}]$ is   a monochromatic triangle, say red, then $(G,\tau)$ contains  a red $W_4=(v_1, v_3, v_2, v_4; v_5)$, where $v_1, v_2\in V_1$, $v_3, v_4\in V_2$ and $v_5\in V_3$. Thus  $\mathcal{R}[\{a_1,a_2,a_3\}]$ is not a monochromatic triangle. Let $A_1,A_2,A_3$ be a permutation of $V_1,V_2,V_3$ such that $A_2$ is, say, blue-complete, to $A_1\cup A_3$ in $G$. Then $A_1$ must be  red-complete to $A_3$. Let $A:=V(G)\less(A_1\cup A_2\cup A_3\cup X)$. Note that there is a red $C_4$ using edges between $A_1$ and $A_3$, and there is a blue $C_4$ using edges between $A_1$ and $A_2$. Since $(G,\tau)$ has neither red nor blue $W_4$, we see that    no vertex in $A$ is red-complete to $A_1\cup A_3$, and no vertex in $A$ is blue-complete to $A_1\cup A_2$ or $A_2\cup A_3$. This implies that $A$ must be red-complete to $A_2$. Note that every vertex in $A$ is   blue-complete to either $A_1$ or   $A_3$.   Let $A^* := \{v\in A \mid  v~\text{is blue-complete to}~A_1  \text{ under  }  \tau\}$. Then   $A\less A^*$ is blue-complete to $A_3$ under  $\tau$.  
Note that $A_1$ is blue-complete to $A^*$ and $A_2$ is red-complete to $A^*$. By Claim~\ref{Two},   $(G[A^*],\tau)$ has neither  red nor blue $P_3$. By Claim \ref{HH}, $|A^* |\le f(k-2)$. Similarly, $|A\less A^*|\le f(k-2)$.  By Claim \ref{V1},  $|A_i|\le |V_1|\le  f(k-2)-|X|$ for all $i\in\{1,2,3\}$. By ($\ast$), $5f(k-2)\le f(k)$ for all $k\ge3$. But then 
\[
|G|=|A^* |+|A\less A^*|+|A_1|+|A_2|+|A_3|+|X|\le 5f(k-2)-2|X|\le f(k),
\]
 contrary to the fact that $|G|= f(k)+1$.\qed\\

By Claim~\ref{V3}, $|V_3|\le 1$. By   Claim \ref{V1}, $|V_1\cup X|\le f(k-2)$ and so $|V_2|\le f(k-2)$. By ($\ast$) again, $5f(k-2)\le f(k)$ for all $k\ge3$. But then 
\[
\begin{split}
|G|&=|V_1\cup X|+|V_2|+|V_3|+\cdots+|V_p|\\
&\le 2f(k-2)+(p-2)\\
&\le 2f(k-2)+12\\
&\le 5f(k-2)\\
&\le f(k),
\end{split}
\]
 because $p\le 14$ and $f(k-2)\ge4$, contrary to the fact that $|G|= f(k)+1$.\medskip
 
This completes the proof of Theorem~\ref{thm}.\qed\\

\section*{Acknowledgements}

 The authors thank both referees for their careful reading and helpful comments.

\end{document}